\documentclass[a4paper,10pt]{article}
\usepackage[utf8]{inputenc}

\usepackage{graphicx}
\usepackage{caption}
\usepackage{subcaption}

\usepackage{blkarray}
\usepackage{amsmath}
\usepackage{amssymb}
\usepackage{amsthm}
\usepackage{accents}

\usepackage{natbib}

\usepackage{color}

\newtheorem*{assumption1}{A1($q$)}
\newtheorem{remark}{Remark}
\newtheorem{theorem}{Theorem}

\numberwithin{equation}{section}

\title{A Numerical Truncation Approximation with\\ A Posteriori  Error Bounds for the Solution of Poisson's Equation}
\author{Saied Mahdian$^{1}$,  Peter W Glynn$^{1}$, Yuanyuan Liu$^{2}$\\ 
$^{1}$Stanford University, USA\\
$^{2}$Central South University, China
}
\date{}

\begin{document}

\maketitle

\begin{abstract}
The solution to Poisson's equation arises in many Markov chain and Markov jump process settings, including that of the central limit theorem, value functions for average reward Markov decision processes, and within the gradient formula for equilibrium Markovian rewards. In this paper, we consider the problem of numerically computing the solution to Poisson's equation when the state space is infinite or very large. In such settings, the state space must be truncated in order to make the problem computationally tractable. In this paper, we provide the first truncation approximation solution to Poisson's equation that comes with provable and computable a posteriori error bounds. Our theory applies to both discrete-time chains and continuous-time jump processes. Through numerical experiments, we show our method can provide highly accurate solutions and tight bounds.
\end{abstract}

\section{Introduction}

Let $X = (X_n: n \geq 0)$ be an irreducible recurrent discrete-time Markov chain taking values in a large or countably infinite discrete state space $S$. If $P = (P(x,y): x,y \in S)$ is the one-step transition matrix, we say that $g$ is a \emph{solution to Poisson's equation for forcing function} $r_c: S \to \mathbb{R}$ if  
\begin{align*}
\sum_{y \in S} P(x,y) |g(y)| < \infty
\end{align*}
for each $x \in S$, and
\begin{align}
(P - I) g = -r_c,
\label{eq:poissoneq}
\end{align}
where we adopt the convention throughout this paper of encoding all functions as column vectors. If $\nu = (\nu(x), x \in S)$ is the unique (up to a multiplicative constant) non-trivial non-negative invariant measure associated with $X$, we assume that
\begin{align}
\sum_{x \in S} \nu(x) |r_c(x)| < \infty
\label{eq:poisson_cond1}
\end{align}
and
\begin{align}
\sum_{x \in S} \nu(x) r_c(x) = 0.
\label{eq:poisson_cond2}
\end{align}
This is the standard solvability condition under which (\ref{eq:poissoneq}) is guaranteed to have a solution, regardless of whether $X$ is positive or null recurrent; see \citet{glynn2023solution}.  
 It is shown there that in the presence of (\ref{eq:poisson_cond1}) and (\ref{eq:poisson_cond2}), that if one fixes $z \in S$ and sets $\tau(z) = \inf\{n \geq 1: X_n = z \}$, then $g^* = (g^*(x): x\in S)$ defines a solution to Poisson's equation for forcing function $r_c$, where 
\begin{align*}
g^*(x) = \mathbb{E}_x \sum_{j=0}^{\tau(z)-1} r_c(X_j)
\end{align*}
for $x \in S$, where $\mathbb{E}_x(\cdot)$  $(\mathbb{P}_x(\cdot))$ represents the expectation (probability), conditional on $X_0 = x \in S$.  

The solution to Poisson's equation arises in many settings. It has long been known that the value function that comes up within average reward optimality for Markov decision processes is the solution to Poisson's equation for the Markov chain associated with the optimal policy; see \citet{ross2014introduction}.  
It also plays a key role within the expression for the equilibrium gradient for a Markov chain; see \citet{rhee2022lyapunov}.  
It furthermore is fundamental to the central limit theorem for additive functionals of Markov chains; see \citet{maigret1978theoreme} and \citet{glynn1996liapunov}. 

In view of the importance of Poisson's equation, the numerical computation of $g^*$ arises as a key computational issue. 
Note that when $X$ has a very large or infinite state space, the state space $S$ must be truncated in order to render numerical computation via solution of linear systems of equations computationally tractable. In such settings, the development of associated error bounds becomes relevant. One approach to computing error bounds involves development of \emph{a priori} bounds. Such error bounds can be computed prior to the computation itself, and are useful when one is choosing the size of the truncation set to be used. This approach is followed in the work of \citet{herve2023computable}; see also the references therein. 

The second approach to developing error bounds is to consider \emph{a posteriori} error bounds. Such bounds use the computed solution of the truncated system, and are only available after the computation has been completed. They tend to be (much) tighter than a priori bounds, and are useful in assessing the errors associated with the computed solution to the truncated system.

The current paper develops, to the best of our knowledge, the first truncation approximation to $g^*$ in terms of computable a posteriori error bounds; the bounds depend on knowledge of a suitable Lyapunov function (See Remarks \ref{rmk5} and \ref{rmk7} for the definitions of our truncation approximations). \citet{liu2022augmented} also proposed a truncation approximation to $g^*$ but without accompanying a posteriori error bounds. Our a posteriori error bounds are provably convergent when the truncation set grows to the entire state space, and our approach is valid for both null recurrent and positive recurrent chains.  
As noted in \citet{infanger2022new}, 
the extension to the null recurrent setting is necessary if one wishes to extend the bounds of this paper to the general Markov jump process setting. 
Our numerical examples show that our method can be highly effective and accurate with moderate truncation set sizes.

In Section \ref{sec:markov_jump_process}, we show how to reduce the computation of the solution to Poisson's equation for jump processes to the discrete time setting, provided that we are willing to cover both the null recurrent and positive recurrent cases.
In Section \ref{sec:bound_expectedreward}, we show how to approximate and bound the expected reward accumulated to the hitting time of a singleton state $z \in S$, when $X$ is a recurrent irreducible discrete-time Markov chain. Section \ref{sec:bound_discretetime_chain} applies the theory of Section \ref{sec:bound_expectedreward} to approximate and bound the solution to Poisson's equation for recurrent irreducible discrete-time Markov chains, while Section \ref{sec:bound_poisson_markov_jump} uses the theory of Section \ref{sec:bound_discretetime_chain} to approximate and bound the solution to Poisson's equation for irreducible positive recurrent Markov jump processes. Section \ref{sec:computational_results} covers computational results. 

\section{Reduction of Poisson's Equation for Markov Jump Processes to the Discrete Time Setting}\label{sec:markov_jump_process}

Let $X = (X(t): t \geq 0)$ be an $S$-valued non-explosive recurrent Markov jump process with rate matrix $Q = (Q(x,y): x,y \in S)$. We say that $h = (h(x): x \in S)$ is a \emph{solution of Poisson's equation with forcing function} $s_c: S \to \mathbb{R}$ if
\begin{align}
\sum_{y \in S} |Q(x,y)| \cdot |h(y)| < \infty
\label{eq:poisson_cont_cond1}
\end{align}
for $x \in S$ and
\begin{align}
Q h = -s_c.
\label{eq:poisson_conteq}
\end{align}
Let $R = (R(x,y): x,y \in S)$ be the one-step transition matrix associated with the embedded discrete time Markov chain for $X$, so that
\begin{align*}
R(x,y) = \begin{cases}
            Q(x,y) /\lambda(x) &, x \neq y \\
            0                  &, x =  y,
         \end{cases}
\end{align*}
where $\lambda(x) \overset{\Delta}{=} -Q(x,x)$ is the jump rate out of state $x \in S$. Note that (\ref{eq:poisson_cont_cond1}) is equivalent to requiring that
\begin{align*}
\sum_{y \in S} R(x,y) |h(y)| < \infty
\end{align*}
for $x \in S$, while (\ref{eq:poisson_conteq}) is equivalent to asserting that 
\begin{align*}
(R - I) h = - s'_c,
\end{align*}
where
\begin{align*}
 s'_c(x) = \frac{s_c(x)}{\lambda(x)}
\end{align*}
for $x \in S$. In other words, $h$ is a solution for Poisson's equation (in discrete time) for forcing function $s'_c$.

The solvability conditions (\ref{eq:poisson_cond1}) and (\ref{eq:poisson_cond2}) translate into the requirement that
\begin{align}
\sum_{x \in S} \frac{\nu(x)}{\lambda(x)} |s_c(x)| < \infty
\label{eq:poisson_cont_cond2}
\end{align}
and
\begin{align}
\sum_{x \in S} \frac{\nu(x)}{\lambda(x)} s_c(x) = 0.
\label{eq:poisson_cont_cond3}
\end{align}
Put $\eta(x) = \nu(x)/\lambda(x)$ for $x \in S$, and observe that (\ref{eq:poisson_cont_cond2}) and (\ref{eq:poisson_cont_cond3}) just assert that
\begin{align}
\sum_{x \in S} \eta(x) |s_c(x)| < \infty
\label{eq:poisson_cont_cond4}
\end{align}
and
\begin{align}
\sum_{x \in S} \eta(x) s_c(x) = 0,
\label{eq:poisson_cont_cond5}
\end{align}
where $\eta = (\eta(x): x \in S)$ is known to be the unique non-negative invariant measure associated with $Q$ that satisfies $\eta Q = 0$ (where we always write mass functions on states as row vectors). 
 The solvability conditions (\ref{eq:poisson_cont_cond4}) and (\ref{eq:poisson_cont_cond5}) are the standard solvability criteria for Poisson's equation associated with jump processes in continuous time; see \citet{glynn2023solution}.  

\section{Bounding Expected Reward Cumulated to Time $\tau(z)$}\label{sec:bound_expectedreward}

To obtain a good bound on $g^*$ that can use a suitable truncation, we will use Lyapunov bounds. Fix a finite subset $K \subseteq S$.

\begin{assumption1}
Given a function $q: S \to \mathbb{R}_+$, we say that condition A1($q$) holds if there exists $v: S \to \mathbb{R}_+$ and $c < \infty$ for which 
\begin{align*}
(P v) (x) \leq v(x) - q(x) + c I(x \in K)
\end{align*}
for $x \in S$.
\end{assumption1}

\begin{remark}
Given A1(q), we say that $v$ is a \emph{Lyapunov function} for bounding cumulative rewards associated with $q$.
\end{remark}

Assume that we have a computational platform capable of solving linear systems of equations corresponding to a finite truncation set $A$ containing $K$. We now write $P$ in block partitioned form, namely
\begin{align*}
P = \begin{blockarray}{cccc} 
      & K & A'& A^c \\
      \begin{block}{c[ccc]}
      K   & P_{11} & P_{12} &P_{13} \\
      A'  & P_{21} & P_{22} &P_{23} \\
      A^c & P_{31} & P_{32} &P_{33} \\
      \end{block}
    \end{blockarray}
\end{align*}
where $A' = A - K$. Similarly, we can write $q$ and $v$ in block partitioned form, namely
\begin{align*}
q = \begin{blockarray}{cc} 
      \begin{block}{c[c]}
      K   & q_{1} \\
      A'  & q_{2} \\
      A^c & q_{3} \\
      \end{block}
    \end{blockarray},
v = \begin{blockarray}{cc} 
      \begin{block}{c[c]}
      K   & v_{1} \\
      A'  & v_{2} \\
      A^c & v_{3} \\
      \end{block}
    \end{blockarray}.
\end{align*}
Select $z \in K$ and let $\psi(z) = \inf \{n \geq 0: X_n = z \}$. Observe that $\psi(z) = \tau(z)$ when $X_0 \neq z$. Put 
\begin{align*}
f(x) = \mathbb{E}_x \sum_{j = 0}^{\psi(z)-1} q(X_j)
\end{align*}
for $x \in A$. 
Of course, $f(z) = 0$. 
Let $T = \inf \{ n \geq 0: X_n \in A^c  \}$ and $T_K = \inf\{ n \geq 1: X_n \in K  \}$. 
Then for $x \in A - \{z\}$, by applying the strong Markov property at time $T_K \leq \tau(z)$, we arrive at the equality
\begin{align}
f(x) &= \mathbb{E}_x \sum_{j=0}^{T_K-1} q(X_j) +  \sum_{y \in (K - \{z\}) } \mathbb{P}_x(X_{T_K} = y) f(y) \nonumber \\
     &= \mathbb{E}_x \sum_{j=0}^{(T_K\wedge T)-1} q(X_j) + \mathbb{E}_x \sum_{j=T}^{T_K-1} q(X_j)I(T<T_K) \nonumber \\ 
     &+ \sum_{y \neq z} \mathbb{P}_x(X_{T_K} = y, T_K < T) f(y) + \sum_{y \neq z} \mathbb{P}_x(X_{T_K} = y, T < T_K) f(y).
\label{eq:f_eqn}
\end{align}

Set $\tilde{K} =  K-\{z\}$, and 
let $G = (G(x,y): x,y \in \tilde{K})$ be the matrix with entries given by
\begin{align*}
G(x,y) = \mathbb{P}_x(X_{T_K} = y, T_K < T).
\end{align*}
Then, it is straightforward to verify that
\begin{align*}
 G(x,y) &= P(x,y) + \left( P_{12} \sum_{n=0}^{\infty}P^n_{22} P_{21}  \right) (x,y)\\
        &= P(x,y) + \left(P_{12}(I - P_{22})^{-1} P_{21}  \right) (x,y),
\end{align*}
where we used the irreducibility of $X$ to ensure that $\sum_{n=0}^{\infty} P^n_{22} < \infty$ (and hence $(I - P_{22})^{-1}$ exists); see \citet{kemeny1960finite}.  
For $x \in A'$, 
\begin{align*}
\mathbb{P}_x(X_{T_K} = y, T_K < T) = \left((I-P_{22})^{-1} P_{21} \right)(x,y).
\end{align*}
Also,  
\begin{align*}
k(x,q) &\overset{\Delta}{=} \mathbb{E}_x \sum_{j=0}^{(T_K \wedge T)-1} q(X_j) =   \sum_{j=0}^{+\infty} \mathbb{E}_x q(X_j) I((T_K \wedge T)-1\geq j) =  \\
       &= \begin{cases} 
            q_1(x) + \left(P_{12}(I - P_{22})^{-1} q_2 \right) (x),  x\in K  \\ 
            \left( (I - P_{22})^{-1} q_2 \right)(x), x \in A'.
          \end{cases}
\end{align*}

In the presence of Assumption A1($q$), it is well known (see, for example, \citet{meyn2009markov}, p.344)  
that the Lyapunov bound 
\begin{align*}
 \mathbb{E}_x \sum_{j=0}^{T_K-1} q(X_j) \leq v(x)
\end{align*}
holds for  $x \in K^c$, and hence
\begin{align*}
k'(x,q) &\overset{\Delta}{=} \mathbb{E}_x \sum_{j=T}^{T_K-1} q(X_j) I(T < T_K) \\
        &\leq \mathbb{E}_x v(X_T)I(T < T_K) \\ 
        & = \begin{cases}
              (P_{13}v_3)(x) + (P_{12}(I - P_{22})^{-1}P_{23} v_3)(x), x\in K\\
              (I - P_{22})^{-1}(P_{23}v_3)(x),  x \in A'
            \end{cases}\\
        &\overset{\Delta}{=}  \tilde{k}'(x,q)
\end{align*}
for $x \in A$.

\begin{remark}
A1($q$) implies that $(P_{13}v_3)(x) < \infty$ for $x \in K$ and $(P_{23}v_3)(x) < \infty$ for $x \in A'$. 
\end{remark}

\begin{remark}
 Note that since $A^c$ can be infinite, we can not assume that $P_{13}v_3$ and $P_{23} v_3$ can be computed in closed form in complete generality. Of course, if each row of $P_{13}$ and $P_{23}$ has finitely many non-zeros, this can be done, or if the Markov chain has special structure it may be possible. Even if $P_{13}v_3$ and $P_{23}v_3$ cannot be computed in closed form, we often can upper bound $P_{13}v_3$ and $P_{23} v_3$, thereby yielding a suitable upper bound $\tilde{k}'(\cdot,q)$ on $k'(\cdot,q)$.
\end{remark}

Hence, it follows from (\ref{eq:f_eqn}) that for $x \in \tilde{K}$,
\begin{equation}
\begin{aligned}
f(x) \leq&  k(x,q) + \tilde{k}'(x,q) + \sum_{y\in \tilde{K}} G(x,y) f(y) \\ 
     &+ \sum_{y\in \tilde{K}} \mathbb{P}_x(X_{T_K}=y, T<T_K) \max_{w \in \tilde{K}} f(w) \\
     \leq&  k(x,q) + \tilde{k}'(x,q) + \sum_{y\in \tilde{K}} G(x,y) f(y) + \mathbb{P}_x( T < T_K) \cdot \| f\|_{\tilde{K}}
\end{aligned}
\label{eq:f_bound1}
\end{equation}
where we adopt the notation $\| \rho \|_{\tilde{K}} \overset{\Delta}{=} \max\{|\rho(x)|: x \in \tilde{K} \}$ for a generic function $\rho: \tilde{K} \to \mathbb{R}$. Set 
\begin{align*} 
\beta(x,q) = k(x,q) + \tilde{k}'(x,q), \xi(x) = 1 - \mathbb{P}_x(X_{T_K}=z, T_K<T) - \sum_{y \in \tilde{K}} G(x,y),
\end{align*}
$f = (f(x): x\in \tilde{K})$, $k(q) = (k(x,q), x \in \tilde{K})$, $\beta(q) = (\beta(x,q): x \in \tilde{K})$ and $\xi = (\xi(x): x \in \tilde{K})$. Then (\ref{eq:f_bound1}) can be viewed as a matrix-vector inequality indexed by states in $\tilde{K}$, namely 
\begin{align*}
 f \leq \beta(q) + G f +  \xi \cdot \| f\|_{\tilde{K}}.
\end{align*}
This leads to the inequality 
\begin{align*}
(I - G) f \leq \beta(q) + \xi \cdot \| f \|_{\tilde{K}}.
\end{align*}
Since $(I - G)^{-1}$ exists and is non-negative, it follows that
\begin{align}
 f \leq (I - G )^{-1} \beta(q) + (I - G)^{-1} \xi \cdot \| f \|_{\tilde{K}}.
\label{eq:f_bound2}
\end{align}
Hence, 
\begin{align*}
 \| f \|_{\tilde{K}} \leq \| (I - G )^{-1} \beta(q)\|_{\tilde{K}}  +  \| (I - G)^{-1} \xi \|_{\tilde{K}} \cdot \| f \|_{\tilde{K}}
\end{align*}
so that when $A$ is large enough that $\| (I - G)^{-1} \xi \|_{\tilde{K}} < 1$, 
\begin{align*}
\| f \|_{\tilde{K}}  \leq \frac{\| (I - G )^{-1} \beta(q)\|_{\tilde{K}}}{ 1 - \| (I - G)^{-1} \xi \|_{\tilde{K}}}.
\end{align*}
We conclude from (\ref{eq:f_bound2}) that 
\begin{align*}
 f \leq& (I - G)^{-1}\beta(q) + (I - G)^{-1} \xi \frac{\|(I-G)^{-1}\beta(q)\|_{\tilde{K}}}{1-\|(I-G)^{-1}\xi\|_{\tilde{K}}} \\ 
  \overset{\Delta}{=}& \tilde{\kappa}(q)
\end{align*}

On the other hand, the non-negativity of $q$ and (\ref{eq:f_eqn}) together imply that
\begin{align}
  f \geq k(q) + G\; f.
\label{eq:f_lowerbound_eqn}
\end{align}
Since $(I-G)^{-1} = \sum_{n=0}^{\infty} G^n \geq 0$, (\ref{eq:f_lowerbound_eqn}) implies that
\begin{align*}
f &\geq (I-G)^{-1}k(q) \\ 
&\overset{\Delta}{=} \undertilde{\kappa}(q).
\end{align*}
We now have upper and lower bounds on $f(x)$ for $x \in \tilde{K}$. On the other hand, for $x \in A'$,
using (\ref{eq:f_eqn}) and 
the same reasoning as above establishes that,
\begin{align*}
 f(x) \leq& k(x,q) + \tilde{k}'(x,q) \\
                  &+ \sum_{y \in \tilde{K}} \mathbb{P}_x(X_{T_K} = y, T_K < T) \tilde{\kappa}(y,q) \\
                  &+ \sum_{y \in \tilde{K}} \mathbb{P}_x(X_{T_K} = y, T < T_K) \tilde{\kappa}(y,q) \\
                  &\leq \beta(x,q) + \sum_{y \in \tilde{K}}\left( (I - P_{22})^{-1} P_{21}  \right) (x,y) \tilde{\kappa}(y,q) + \xi(x)\cdot \| \tilde{\kappa}(q) \|_{\tilde{K}}
\end{align*}
where we have extended $\xi(\cdot)$ to $A'$ with the definition
\begin{align*} 
\xi(x) = \mathbb{P}_x(T < T_K) = 1 - \sum_{y \in K} \left( (I - P_{22})^{-1} P_{21} \right)(x,y)
\end{align*}
for $x \in A'$. It follows that for $x \in A'$,
\begin{align*}
f(x) \leq \tilde{\kappa}(x,q)
\end{align*}
where $\tilde{\kappa}(\cdot,q)$ is extended to $A'$ via
\begin{align*} 
 \tilde{\kappa}(x,q) = \beta(x,q) + \sum_{y \in \tilde{K}} \left( (I - P_{22})^{-1} P_{21} \right)(x,y) \tilde{\kappa}(y,q) + \xi(x) \| \tilde{\kappa}(q) \|_{\tilde{K}}.
\end{align*}

For the lower bound on $f(x)$ when $x \in A'$, the same argument as for $f(x)$ with $x \in \tilde{K}$ shows that for $x \in A'$,
\begin{align*}
f(x) \geq \undertilde{\kappa}(x,q)
\end{align*}
where
\begin{align*}
 \undertilde{\kappa}(x,q) \overset{\Delta}{=} k(x,q) + \sum_{y \in \tilde{K}} \left( (I - P_{22})^{-1} P_{21}\right) (x,y) \undertilde{\kappa}(y,q)
\end{align*}
for $x \in A'$.

We have therefore arrived at the following theorem.

\begin{theorem}
Suppose that $X$ is an irreducible recurrent Markov chain satisfying A1($q$) for some function $q: S \to \mathbb{R}_+$. Then, 
\begin{align*}
\undertilde{\kappa}(x,q) \leq \mathbb{E}_x \sum_{j=0}^{\tau(z)-1} q(X_j) \leq \tilde{\kappa}(x,q)
\end{align*}
for each $x \in A - \{z\}$, provided A is large enough that $\|(I - G)^{-1}\xi \|_{\tilde{K}} < 1$. 
\label{thm:expectedreward_bound}
\end{theorem}

Since both $\undertilde{\kappa}(\cdot,q)$ and $\tilde{\kappa}(\cdot,q)$ are numerically computable, this provides us with usable bounds on $f(x)$ for $x \in A$.

\begin{remark}
 We note that A1($q$) does not imply that $X$ is positive recurrent. For example, when $S = \mathbb{Z}_+$ and $q(x) \to 0$ as $x \to \infty$, one could have
\begin{align}
\mathbb{E}_x \sum_{j=0}^{\tau(z)-1} q(X_j) < \infty
\label{eq:f_cond1}
\end{align} 
for each $x \in S$, without $\mathbb{E}_z \tau(z) < \infty$. We also note that (\ref{eq:f_cond1}) holds for $x \in S$ if and only if A1($q$) is valid. (We may take $v(x) = \mathbb{E}_x \sum_{j=0}^{\tau(z)-1} q(X_j). $)
\end{remark}

\begin{remark}
Consider a sequence $(A_n: n \geq 1)$ of truncation sets for which $K \subset A_1 \subset A_2 \subset \cdots$ with $A_n \nearrow S$. Write $\tilde{\kappa}_n(x,q)$ and $\undertilde{\kappa}_n(x,q)$ to emphasize the dependence of the upper and lower bounds on the choice of the truncation set $A_n$, and put $T_n = \inf\{ m \geq 0: X_m \in A_n^c \}$ for $n \geq 1$. Then, 
\begin{align*}
\undertilde{\kappa}_n(q) = (I-G_n)^{-1}k_n(q) = \sum_{j=0}^{\infty} G_n^j \; k_n(q)\leq f. 
\end{align*}
On the other hand, the first line of (\ref{eq:f_eqn}) implies that
\begin{align*}
 f = (I - H)^{-1} h(q)  = \sum_{j=0}^{\infty} H^j\;h(q),
\end{align*}
where $H=(H(x,y): x,y \in \tilde{K})$ has entries given by $H(x,y) = \mathbb{P}_x(X_{T_K}=y)$, and $(h(q))(x) = \mathbb{E}_x \sum_{j=0}^{T_K-1} q(X_j)$.
Note that because $T_n \nearrow \infty$ as  $n \to \infty$, it follows that
\begin{align*}
k_n(x,q)  \overset{\Delta}{=} \mathbb{E}_x \sum_{j=0}^{(T_K \wedge T_n)-1} q(X_j) \nearrow \mathbb{E}_x \sum_{j=0}^{T_K-1} q(X_j) =  (h(q))(x)
\end{align*}
via the Monotone Convergence Theorem and $G_n \nearrow H$ as $n \to \infty$.
So, $\undertilde{\kappa}_n(q)$ converges monotonically to $f$ as $n \to \infty$.
We may therefore use $\undertilde{\kappa}_n(x,q)$ as \emph{our truncation approximation} to $f(x)$, since we are then guaranteed that the approximation converges to $f(x)$.

On the other hand, $\xi_n(x) = \mathbb{P}_x (T_n < T_K) \to 0$  as $n \to \infty$. Also, 
$k_n(x,q) \to \mathbb{E}_x \sum_{j=0}^{T_K-1} q(X_j)$
as $n \to \infty$ and
\begin{align*}
\tilde{k}'_n(x,q) \overset{\Delta}{=} \mathbb{E}_x v(X_{T_n}) I(T_n < T_K) \leq \mathbb{E}_x \sum_{j=0}^{\tau(z)-1} v(X_j).
\end{align*}
If condition A1($v$) holds, then $\sum_{j=0}^{\tau(z)-1} v(X_j)$ is $\mathbb{P}_x$-integrable, so $\mathbb{E}_x v(X_{T_n}) I(T_n < T_K) \to 0$ via the Dominated Convergence Theorem. Hence, our upper bound $\tilde{\kappa}_n(x,q)$  converges to $f(x)$ as $n \to \infty$, establishing that our upper bound is convergent under the extra condition A1($v$). 
\label{rmk5}
\end{remark}

\section{Numerical Bounds for Poisson's Equation for Discrete-time Markov Chains}\label{sec:bound_discretetime_chain}

Suppose $r: S \to \mathbb{R}_+$ is a non-negative reward function. We assume that $X$ is irreducible and positive recurrent.
Then, we can normalize the invariant probabilities $\nu$ so that it is a probability $\pi = (\pi(x): x \in S)$.

The solvability conditions (\ref{eq:poisson_cond1}) and (\ref{eq:poisson_cond2}) imply that $\sum_{x\in S} \pi(x) r(x) < \infty$, and that we work with the ``centered reward'' function $r_c: S \to \mathbb{R}$ given by
\begin{align*}
r_c(x) = r(x) - \alpha,
\end{align*}
where $\alpha = \pi r ( = \sum_{y \in S}\pi(y)r(y))$. We now wish to use the machinery of Section \ref{sec:bound_expectedreward} to bound
\begin{align*} 
g^*(x) = \mathbb{E}_x \sum_{j=0}^{\tau(z)-1}r_c(X_j)
\end{align*}
over $x \in A$. Of course we can write
\begin{align*}
g^*(x) = \mathbb{E}_x \sum_{j=0}^{\tau(z)-1}r(X_j) - \alpha\; \mathbb{E}_x \sum_{j=0}^{\tau(z)-1}e(X_j),
\end{align*}
where $e: S \to \mathbb{R}_+$ is given by $e(x) = 1$ for $x \in S$. Note that both $r$ and $e$ are non-negative, so the bounds of Section \ref{sec:bound_expectedreward} apply if A1($r$) and A1($e$) are assumed to hold. Hence, when $\| (I - G)^{-1}\xi \|_{\tilde{K}} < 1$, Theorem \ref{thm:expectedreward_bound} applies, so
\begin{align*}
\undertilde{\kappa}(x,r) \leq \mathbb{E}_x \sum_{j=0}^{\tau(z)-1} r(X_j) \leq \tilde{\kappa}(x,r)
\end{align*}
and
\begin{align*}
\undertilde{\kappa}(x,e) \leq \mathbb{E}_x \tau(z)  \leq \tilde{\kappa}(x,e)
\end{align*}
for $x \in A - \{z\}$. Furthermore, 
\begin{align*}
\mathbb{E}_z \sum_{j=0}^{\tau(z)-1} r(X_j) = r(z) + \sum_{y \neq z} P(z,y)  \mathbb{E}_y \sum_{j=0}^{\tau(z)-1} r(X_j).
\end{align*}
As a result,
\begin{align*}
\undertilde{\kappa}(z,r) = r(z) + \sum_{y \in A, y \neq z} P(z,y) \undertilde{\kappa}(y,r) \leq \mathbb{E}_z \sum_{j=0}^{\tau(z)-1} r(X_j)
\end{align*}
and
\begin{align*}
\tilde{\kappa}(z,r) &= r(z) + \sum_{y \in A, y \neq z} P(z,y) \tilde{\kappa}(y,r)  \\
                    &+ \sum_{y \in A^c} P(z,y) \left( v_3(y) + \| \tilde{\kappa}(r)  \|_{\tilde{K}}   \right) \\ 
                    & \geq \mathbb{E}_z \sum_{j=0}^{\tau(z)-1} r(X_j).
\end{align*}
Similar upper and lower bounds $\tilde{\kappa}(z,e)$ and $\undertilde{\kappa}(z,e)$ for $\mathbb{E}_z \tau(z)$ can be derived, using the Lyapunov function that goes with A1($e$). The regenerative structure of $X$ implies that
\begin{align*}
\alpha = \frac{\mathbb{E}_z \sum_{j=0}^{\tau(z)-1} r(X_j)}{ \mathbb{E}_z \tau(z)};
\end{align*}
see, for example, \citet{smith1955regenerative}. Hence, we obtain the following bounds on $\alpha$:
\begin{align}
\undertilde{\alpha} \overset{\Delta}{=} \frac{\undertilde{\kappa}(z,r)}{\tilde{\kappa}(z,e)} \leq \alpha \leq  \frac{\tilde{\kappa}(z,r)}{\undertilde{\kappa}(z,e)} \overset{\Delta}{=} \tilde{\alpha}
\label{eq:bounda}
\end{align}
Of course, there are many other upper and lower numerically computable bounds for $\alpha$ that have now been derived in the literature; see, for example, \citet{infanger2022new} and the references therein. In what follows, we may use the bounds associated with (\ref{eq:bounda}), or we may use any other upper and lower bounds for $\alpha$ in their place. We therefore arrive at the following computable upper and lower bounds on $g^*$, namely
\begin{align}
\undertilde{\kappa}(x,r) - \tilde{\alpha} \tilde{\kappa}(x,e) \leq g^*(x) \leq \tilde{\kappa}(x,r) - \undertilde{\alpha} \undertilde{\kappa}(x,e)
\label{eq:boundg}
\end{align}
for $x \in A$, with $g^*(z)=0$. The inequalities (\ref{eq:boundg}) represent our bounds on the solution $g^*$ to Poisson's equation.

\begin{remark}
If $r: S \to \mathbb{R}$ is of mixed sign, the easiest way to proceed is to separately derive bounds on the the solutions of Poisson's equation associated with the positive and negative parts of $r$.
\end{remark}

\begin{remark}
To approximate $g^*(x)$, any value in the closed interval 
\begin{align*}
\left[\undertilde{\kappa}(x,r) - \tilde{\alpha} \tilde{\kappa}(x,e), \tilde{\kappa}(x,r) - \undertilde{\alpha} \undertilde{\kappa}(x,e)\right]
\end{align*}
can be used.
However, we may choose  
\begin{align*}
\undertilde{\kappa}(x,r) - \frac{\undertilde{\kappa}(z,r)}{\undertilde{\kappa}(z,e)} \undertilde{\kappa}(x,e)
\end{align*}
as \emph{our truncation approximation} to $g^*(x)$ for $x \in A$, as discussed in Remark $\ref{rmk5}$, because it is constructed from lower bounds that are guaranteed to be convergent under the conditions stated in Remark $\ref{rmk5}$. 
\label{rmk7}
\end{remark}

\section{Numerical Bounds for the Solution of Poisson's Equation for Markov Jump Processes}\label{sec:bound_poisson_markov_jump}

Suppose that $s: S \to \mathbb{R}_+$ is a non-negative reward function. We assume that $(X(t): t \geq 0)$ is an irreducible non-explosive positive recurrent Markov jump process. Then, the invariant measure $\eta = (\eta(x): x \in S)$ has finite total mass, and the solvability conditions for Poisson's equation require that $\eta s = \sum_{x \in S} \eta(x) s(x) < \infty$ and that we work with the centered reward function
\begin{align*}
s_c(x) = s(x) - \delta,
\end{align*} 
where
\begin{align*}
\delta = \frac{\sum_{y \in S} \eta(s) s(y)}{\sum_{y \in S} \eta(y) e(y)} = \frac{\eta s}{\eta e}.
\end{align*}
The centering requirement demands that
\begin{align*}
0 = \sum_{x \in S} \eta(x) [s(x) - \delta] = \sum_{x \in S} \frac{\nu(x)}{\lambda(x)} [s(x) - \delta],  
\end{align*}
where $\nu=(\nu(x): x \in S)$ is the non-negative invariant measure of the embedded discrete-time Markov chain $Y= (Y_n: n \geq 0)$ (having transition matrix R) associated with $(X(t): t \geq 0)$. As discussed in Section \ref{sec:markov_jump_process}, the solution to $Q h = -s_c$ is then given by
\begin{align*}
h^*(x) = \mathbb{E}_x \sum_{j=0}^{\tau(z)-1} \frac{s(Y_j)}{\lambda(Y_j)} - \delta \mathbb{E}_x \sum_{j=0}^{\tau(z)-1} \frac{e(Y_j)}{\lambda(Y_j)} 
\end{align*}
for $x \neq z$, with $h^*(z) = 0$. Put
\begin{align*}
& s'(x) = s(x) / \lambda(x), \\
& e'(x) = e(x) / \lambda(x),
\end{align*}
for $x \in S$. In the presence of A1($s'$) and A1($e'$), Section \ref{sec:bound_expectedreward}'s analysis provides upper and lower bounds such that
\begin{align*}
\undertilde{\kappa}(x,s') \leq \mathbb{E}_x \sum_{j=0}^{\tau(z) -1} s'(Y_j) \leq \tilde{\kappa} (x,s')
\end{align*}
and
\begin{align*}
\undertilde{\kappa}(x,e') \leq \mathbb{E}_x \sum_{j=0}^{\tau(z) -1} e'(Y_j) \leq \tilde{\kappa} (x,e')
\end{align*}
for $x \in A - \{z\}$. As for $\delta$, the principle of regeneration implies that when $(X(t): t\geq 0)$ is positive recurrent, then
\begin{align*}
\mathbb{E}_z \sum_{j=0}^{\tau(z)-1} e'(Y_j) < \infty,  
\end{align*}
and
\begin{align*}
 \delta = \frac{\mathbb{E}_z \sum_{j=0}^{\tau(z)-1} s'(Y_j)}{\mathbb{E}_z \sum_{j=0}^{\tau(z)-1} e'(Y_j)}.
\end{align*}
In Section \ref{sec:bound_discretetime_chain}, we provided upper and lower bounds $\tilde{\kappa}(z,s')$, $\tilde{\kappa}(z,e')$, $\undertilde{\kappa}(z,s')$ and $\undertilde{\kappa}(z,e')$ that are valid so long as $Y$ is recurrent and A1($s'$) and A1($e'$) hold. This leads to the desired bounds on the solution $h^*$ to Poisson's equation for the jump process, namely 
\begin{align}
\undertilde{\kappa}(x,s') - \frac{\tilde{\kappa}(z,s')}{\undertilde{\kappa}(z,e')} \tilde{\kappa}(x,e') \leq h^*(x) \leq \tilde{\kappa}(x,s') - \frac{\undertilde{\kappa}(z,s')}{\tilde{\kappa}(z,e')} \undertilde{\kappa}(x,e')
\label{eq:boundh}
\end{align}
for $x \in A - \{z\}$. (Of course, $h^*(z)=0$.)

\section{Computational Results}\label{sec:computational_results}

In this section, we provide numerical results to illustrate the performance of our bounds. 

\subsection{A Discrete-time Queue}\label{sec:discretetime_queue}
We consider a slotted time queueing model in which the number-in-system process $X = (X_n: n \ge 0)$ evolves according to the recursion:
\begin{align}
  X_{n+1} = \left[ X_n + B_n - D_{n+1} \right]^+
\end{align}
for $n \ge 0$. Here, $B_n$ is the number of customers arriving just after the start of slot $n$, $D_{n+1}$ is the (independent) number of customers served within slot $n$ (departing just prior to the start of slot $n+1$), and $X_n$ is the number of customers in the system at the start of slot $n$. The notation $[y]^+$ means max$\{y,0\}$. 

Suppose that $B_1$ is uniform on $1, 2, 3$ and assume that $D_1$ is geometric with parameter $1- q$ $(0 < q < 1)$. Specifically,
\begin{align}
 \mathbb{P}(D_1 = k) = (1-q) q^{k-1}, \quad k=1,2,\ldots
\end{align}
The stability condition for the queue is then given by $(1-q) \mathbb{E}B_1 = 2(1-q) < 1$. Define $Z_n = B_n - D_{n+1}$. A key property of this model is that for $k \ge 1$,
\begin{align*}
  \mathbb{P}(Z_n = -k) &= (1-q) q^{k-1} \sum_{j=0}^{\infty} P(B_1 = j) q^j \\
  &= (1-q) q^{k-1} (1/3) (q + q^2 + q^3),
\end{align*}
so that $Z_n$'s left tail is geometric. This geometric left tail property is the discrete state space analog to the exponential left tail property that arises in the continuous state space M/G/1 waiting time sequence analysis of \citet{glynn1994poisson}. As in that paper, the moments of $X_{\infty}$ can be explicitly computed via an identical argument, as can the solution to Poisson's equation for reward functions of the form $r(x) = x^m$ for $m \ge 1$. In particular, when $q=0.6$ and $r(x) = x$, E$r(X_{\infty}) = 8/3$ and the solution to Poisson's equation for $r_c$ that vanishes at the origin is given by $g^*(x) = x^2 + 4x$. 

To test our bounds, we use the Lyapunov function $v(x) = 2x^2$ to satisfy both A1($r$) and A1($e$). The set $K$ associated with $v$ and A1($r$) can be taken to be $K = \{0,\ldots, 9\}$  
while the constant $c$ appearing in A1($r$) can be taken as $c = 24.456$.
These choices of $K$ and $c$ rest on using the inequality
\begin{align*}
  (Pv)(x) =& 2 \mathbb{E}(([x+Z_1]^+)^2) \leq 2 \mathbb{E} (x+Z_1)^2 = \\
           & v(x) - r(x) + (4 \mathbb{E}Z_1 + 1) x + 2 \mathbb{E}Z_1^2.
\end{align*}
We choose $c = 2 \mathbb{E}Z_1^2$ and $K^c$ the set of all $x \in S$ such that $(4 \mathbb{E}Z_1 + 1) x + 2 \mathbb{E}Z_1^2 < 0$. 

Since $r(x) \geq e(x)$ for $x \geq 1$ and $(Pv)(0) \leq v(0) + c - 1$, it follows that $v$ also satisfies A1($e$) for $K$. Finally, put $z = 0$.

\subsection{A Network of Queues 1}\label{sec:network_queues1}
Consider a network consisting of two independent M/M/1 queues, in continuous time. (This example is intended to test our bounds in a multi-dimensional setting in which the solution to Poisson's equation for the network can be explicitly computed.) For $i \in \{1,2\}$, $\lambda_i$ is the arrival rate to station $i$, and $\mu_i$ is the station's service rate parameter. Customers completing service at station $i$ immediately leave the network upon completing service. Put $\lambda_1 = 2$, $\mu_1 = 5$, $\lambda_2 = 1$ and $\mu_2 = 3$, and  define $\rho_i = \lambda_i / \mu_i$.  Let $X = (X(t) = (X_1(t), X_2(t)): t \geq 0)$ be the Markov jump process presenting the number of individuals at each station. The state space of $X$ is $S = \{(i_1,i_2): i_1,i_2 \in \mathbb{Z}_+\}$ where $i_j$ is the number of individuals at station $j$.

We consider the Poisson equation $Q h = -s_c$ where $Q$ is the transition rate matrix for $X$, $s(x) = x_1 + x_2$ for $x = (x_1,x_2)\in S$ and $s_c(x) = s(x) - (\gamma_1 + \gamma_2)$ where $\gamma_i = \mathbb{E}X_i(\infty) = \rho_i(1 - \rho_i)^{-1}$. The solution $h^*$ of Poisson's equation $Qh = -s_c$ vanishing at the origin is given by $h^*(x) = h_1(x_1) + h_2(x_2)$, where 
\begin{align*}
 h_i(y) = \frac{y^2 + y}{2(\mu_i - \lambda_i)}. 
\end{align*}
For $s(x)$, we consider the Lyapunov function $v(x)= x_1^2 + x_2^2$. If
\begin{align*}
K = \{(a_1,a_2): (a_1,a_2)\in S,  5 a_1 + 3 a_2   \leq 11 \}, 
\end{align*}
then 
\begin{align*}
 (Q v)(x) \leq -s(x) + 11 I(x \in K)
\end{align*}
for $x \in S$.

Note that $s(x) \geq e(x)$ for $x \neq (0,0)$ and $(Q\;v)(0,0) = 3$, evidently $(Q\;v)(x) \leq -e(x) + 11 I(x \in K)$ for $x \in S$. So, $v$ satisfies A1($e$) for $K$. Finally, we put $z=(0,0)$.

\subsection{A Network of Queues 2}\label{sec:network_queues2}

In this section, we focus on a example in which the solution of Poisson's equation in closed form is unknown. Consider a Jackson network consisting of two single server first-come-first-serve queues. The routing matrix is 
\begin{align*}
 B = \begin{bmatrix}
       \frac{1}{3} & \frac{1}{6} \\
       \frac{1}{3} & \frac{1}{2} 
     \end{bmatrix}.
\end{align*} 
The exogenous arrival rates are $\lambda_1 = 0.75, \lambda_2 = 1$. The service rates are $\mu_1 = \mu_2 = 4$. Let the Markov process $X = (X_t: t \geq 0)$ represent the number of customers at each station. The state space is $S = \{(i_1,i_2): i_1, i_2 \in \mathbb{N}\cup \{0\} \}$.

We consider the Poisson equation $Q h = -s_c$ where again $Q$ is the transition rate matrix and $s(x) = x_1 + x_2$ for $x = (x_1,x_2) \in S$. We have
\begin{align*}
  s_c(x) = s(x) - \frac{\rho_1}{1 - \rho_1} - \frac{\rho_2}{1 - \rho_2}
\end{align*} 
where $\rho_i = \nu_i / \mu_i$ and $\nu_1, \nu_2$ are the solutions to the following system of equations.
\begin{align*}
  \nu  = \lambda + B^T \nu
\end{align*}
We consider the Lyapunov function $v(x) = 5 x_1^2 + 5 x_2^2$ and set 
\begin{align*}
 K = \{ (a_1,a_2)  :  (a_1, a_2) \in S, 4.8333 a_1 + 2.3333 a_2 \leq 42.0833   \}.
\end{align*}
These choices for $v,K$ ensure $A1(s)$ and $A1(e)$ both hold. We also set $z = (0,0)$.

\subsection{Numerical Results}

\begin{figure}[htbp]
\centering
\begin{subfigure}[b]{0.48\textwidth}
\includegraphics[width=\textwidth]{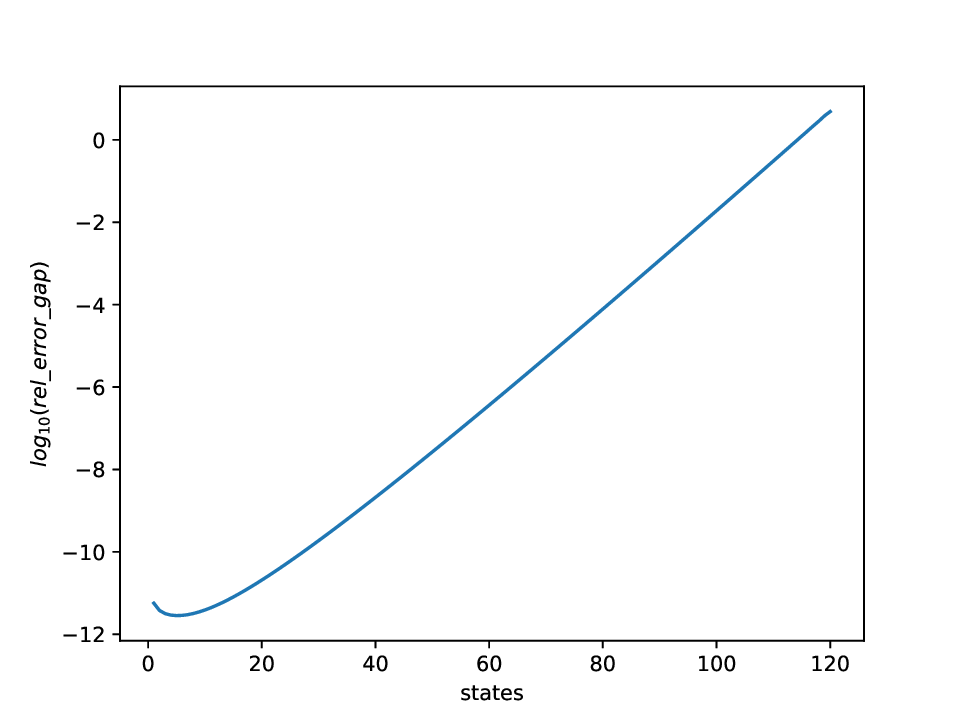}
\caption{Discrete-time queue}
\label{fig1:figa}
\end{subfigure}
\begin{subfigure}[b]{0.48\textwidth}
\includegraphics[width=\textwidth]{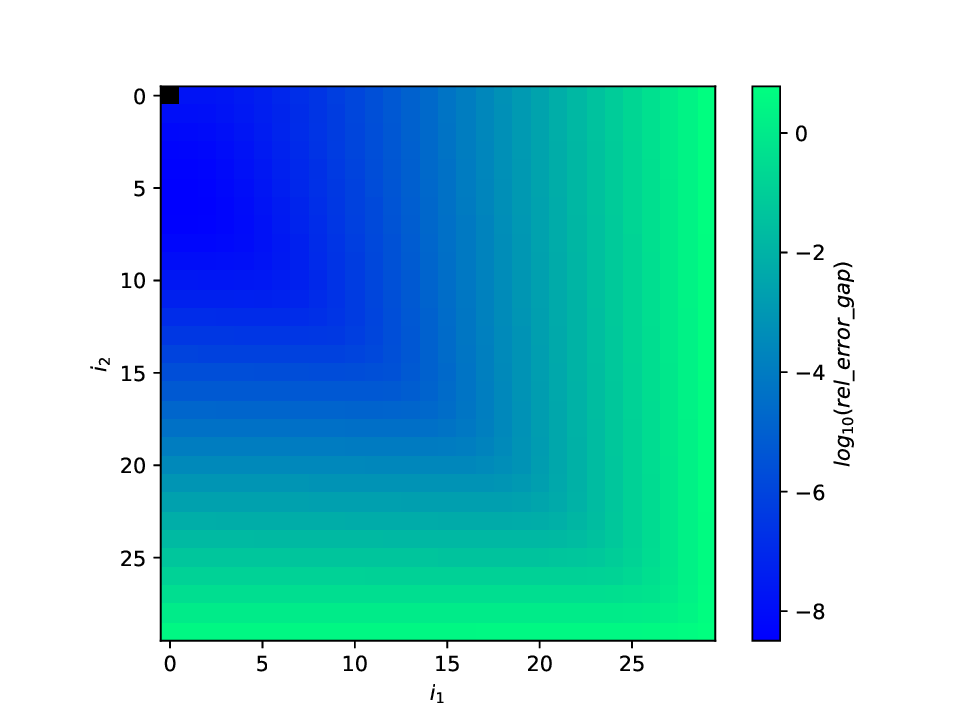}
\caption{Network of queues 1}
\label{fig1:figb}
\end{subfigure}
\caption{Relative error gap for all states in a specific truncation set.}
\label{fig::fig1}
\end{figure}

In this section, we provide numerical results for the settings described in the previous subsections.
Figure \ref{fig::fig1} shows $log_{10}$ of the relative error gap for all states of the truncation set for the settings of Sections \ref{sec:discretetime_queue} and \ref{sec:network_queues1}. The relative error gap for any state $x \in S$ is defined as 
\begin{align*}
rel\_error\_gap(x) = \frac{|upper\_b(x) - lower\_b(x)|}{|exact(x)|}
\end{align*}
when $exact(x) \neq 0$, where $upper\_b,lower\_b,exact$ are the upper bound, lower bound and exact value for the solution to Poisson's equation based on (\ref{eq:boundg}) or (\ref{eq:boundh}) (depending on the setting).
For Figure \ref{fig1:figa}, we set the truncation set to $\{1,2,\cdots,120\}$ and for Figure \ref{fig1:figb}, we set the truncation set to $\{(i_1,i_2) : i_1,i_2 \in S, i_1\leq 30, i_2 \leq 30\}$.
The figures indicate that close to the boundary between $A$ and $A^c$, the relative error gap for the states in $A$ increases and it can be large.  

\begin{figure}[htbp]
\centering
\includegraphics[scale=0.4]{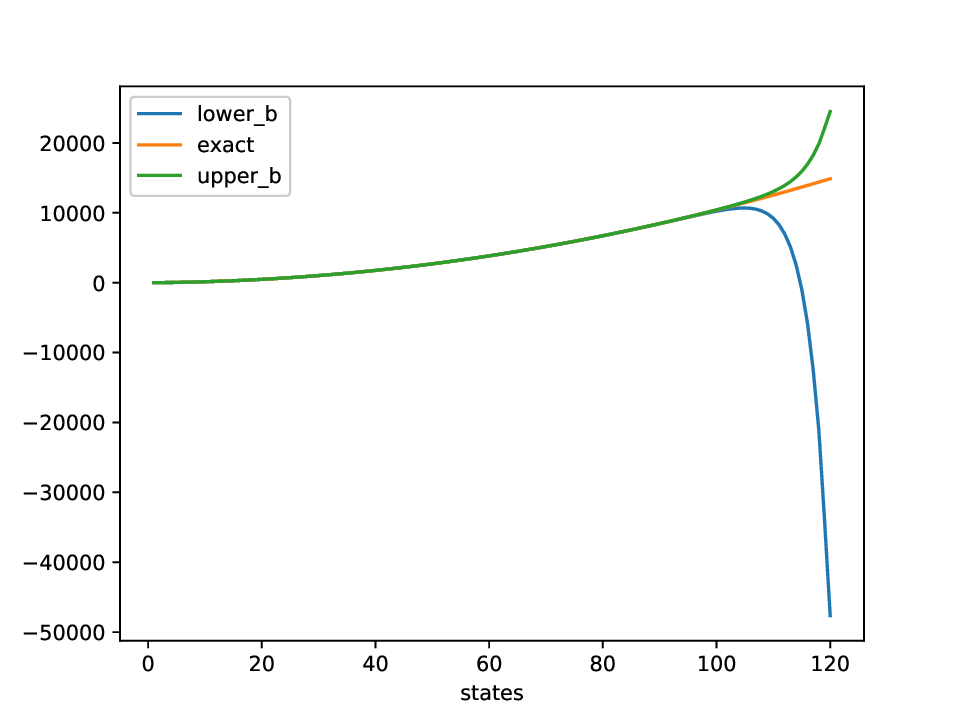}
\caption{Discrete-time queue, value of $upper(\cdot), exact(\cdot), lower(\cdot)$ for each state of the truncation set of Figure \ref{fig1:figa}.}
\label{fig3}
\end{figure}

Figure \ref{fig3} displays $upper\_b(\cdot), lower\_b(\cdot), exact(\cdot)$ for all states of the truncation set (same as in Figure \ref{fig1:figa}) for the discrete-time queue setting. 
Figure \ref{fig3} suggests that close to the boundary between $A$ and $A^c$, the approximation errors of both $upper\_b(\cdot)$ and $lower\_b(\cdot)$ are high. However, the approximation error of $lower\_b(\cdot)$ is much larger as it becomes negative close to the boundary between $A$ and $A^c$. Similar observations were made within the setting of Section \ref{sec:network_queues1}.

\begin{figure}[htbp]
\centering
\begin{subfigure}[b]{0.48\textwidth}
\includegraphics[width=\textwidth]{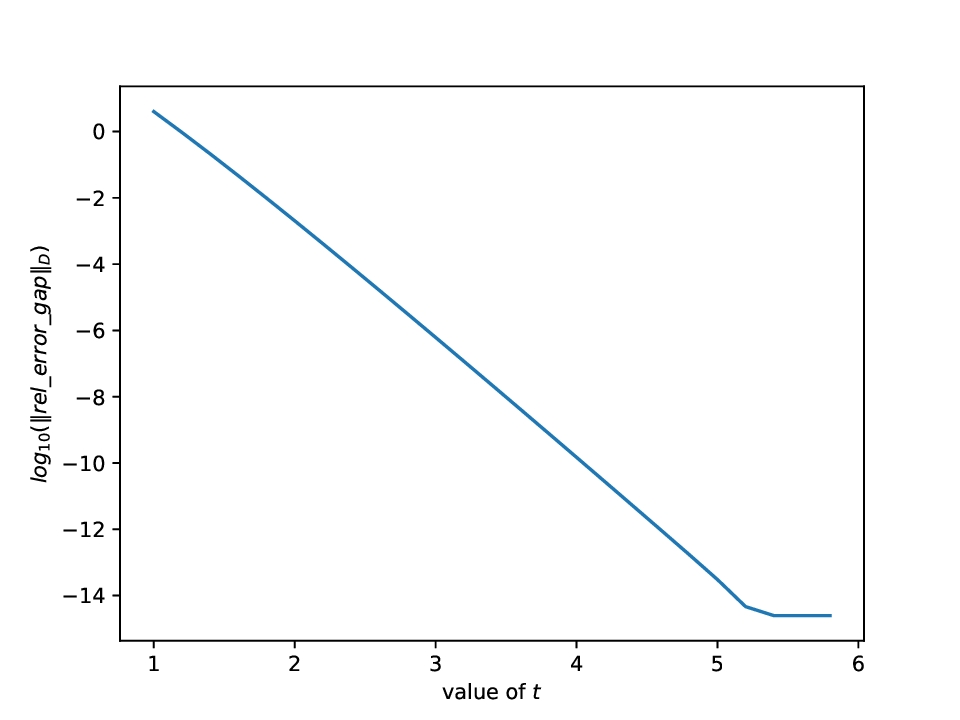}
\caption{Discrete-time queue}
\label{fig2:figa}
\end{subfigure}
\begin{subfigure}[b]{0.48\textwidth}
\includegraphics[width=\textwidth]{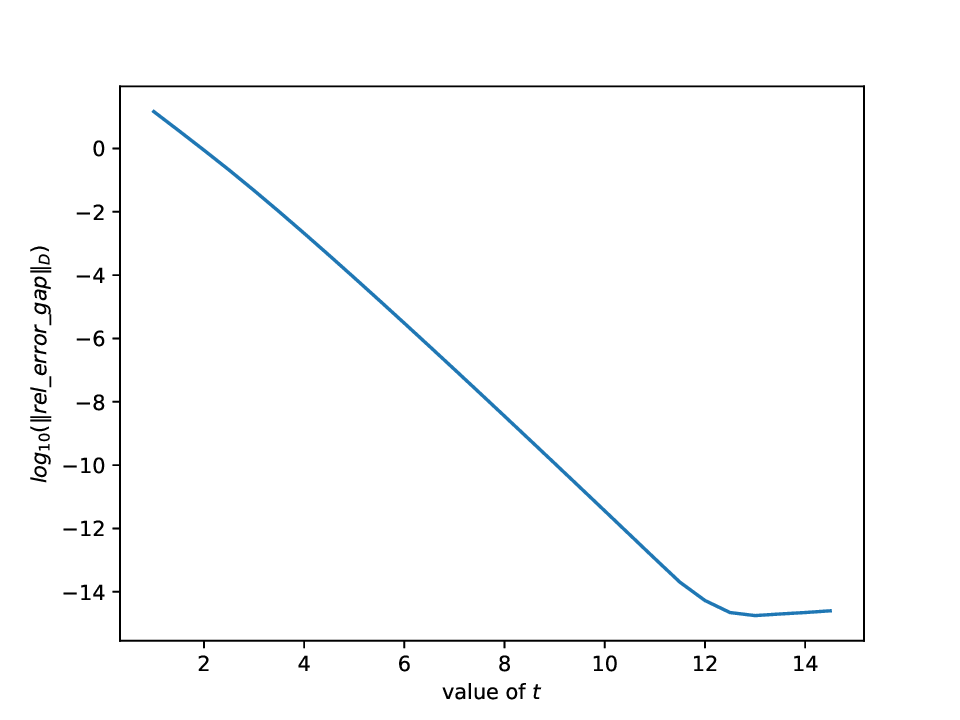}
\caption{Network of queues 1}
\label{fig2:figb}
\end{subfigure}
\caption{Relative error gap versus $t$.}
\label{fig::fig2}
\end{figure}

Next, suppose we are interested in estimating the value of the solution to Poissons's equation over a finite set of states we call $D$, as the size of the truncation set $A_t$ (as a function of $t$) gets large.
Figure \ref{fig::fig2} shows $log_{10}$ of $\|rel\_error\_gap\|_{D} = \max_{x \in D} rel\_error\_gap(x)$ versus $t$ for the settings of Sections \ref{sec:discretetime_queue} and \ref{sec:network_queues1}, where $D = A_1 - \{z\}$, 
and $z$ was specified in Sections \ref{sec:discretetime_queue} and \ref{sec:network_queues1}. 
We do not include $z$ in $D$ because $exact(z) = 0$.
For Figure \ref{fig2:figa}, $A_t = \{x \in S: x \leq \lceil 30 t \rceil \}$ and for Figure \ref{fig2:figb}, $A_t = \{(i_1,i_2) \in S: i_1\leq \lceil 4t \rceil, i_2\leq \lceil 4t \rceil \}$ where $\lceil x \rceil$ is the smallest integer value that is greater or equal to $x$. 
The figures suggest that in both settings $\|rel\_error\_gap\|_{D}$ decreases as a function of $t$ and it can be very small (indicating high accuracy) at moderate values of $t$.

\begin{figure}[htbp]
\centering
\begin{subfigure}[b]{0.48\textwidth}
\includegraphics[width=\textwidth]{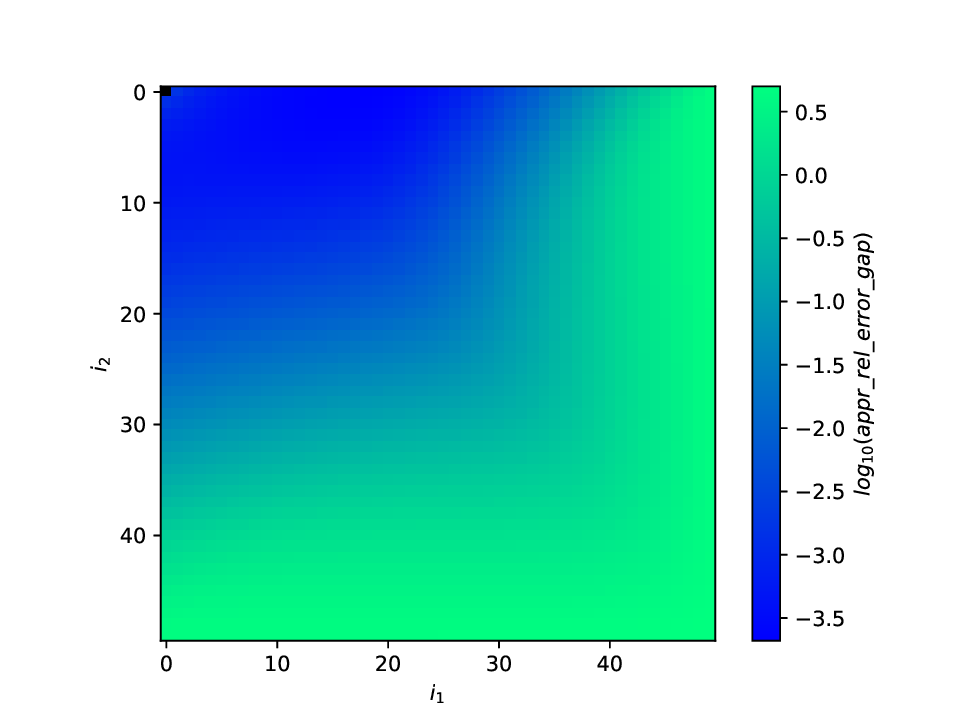}
\caption{Approximate relative error gap}
\label{fig4:figa}
\end{subfigure}
\begin{subfigure}[b]{0.48\textwidth}
\includegraphics[width=\textwidth]{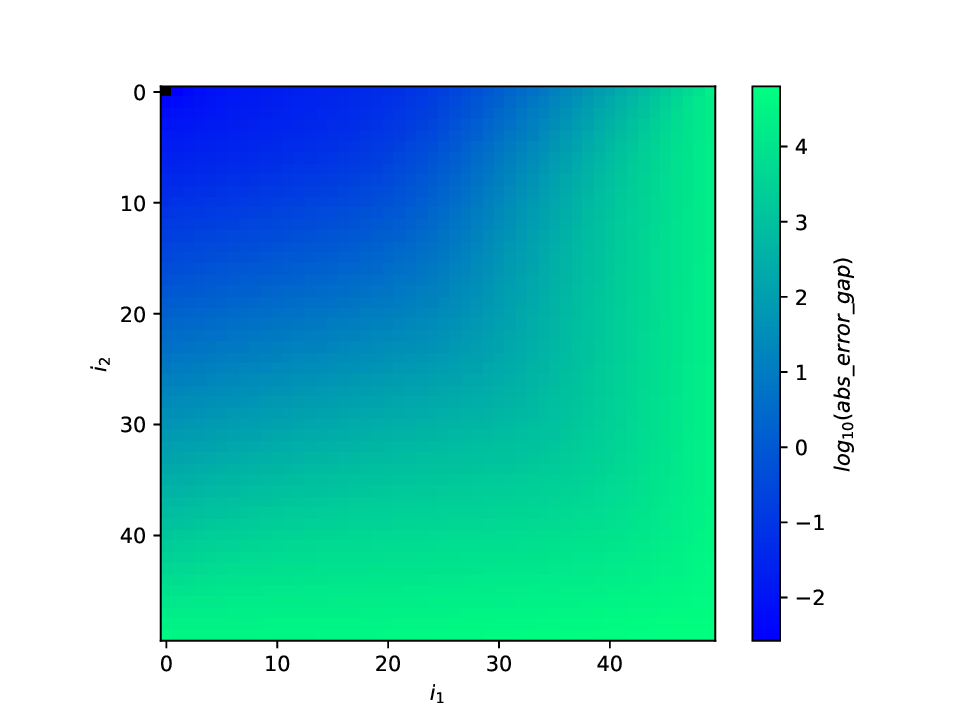}
\caption{Absolute error gap}
\label{fig4:figb}
\end{subfigure}
\caption{Network of queues 2, accuracy measures for a specific truncation set.}
\label{fig::fig4}
\end{figure}

\begin{figure}[htbp]
\centering
\includegraphics[scale=0.4]{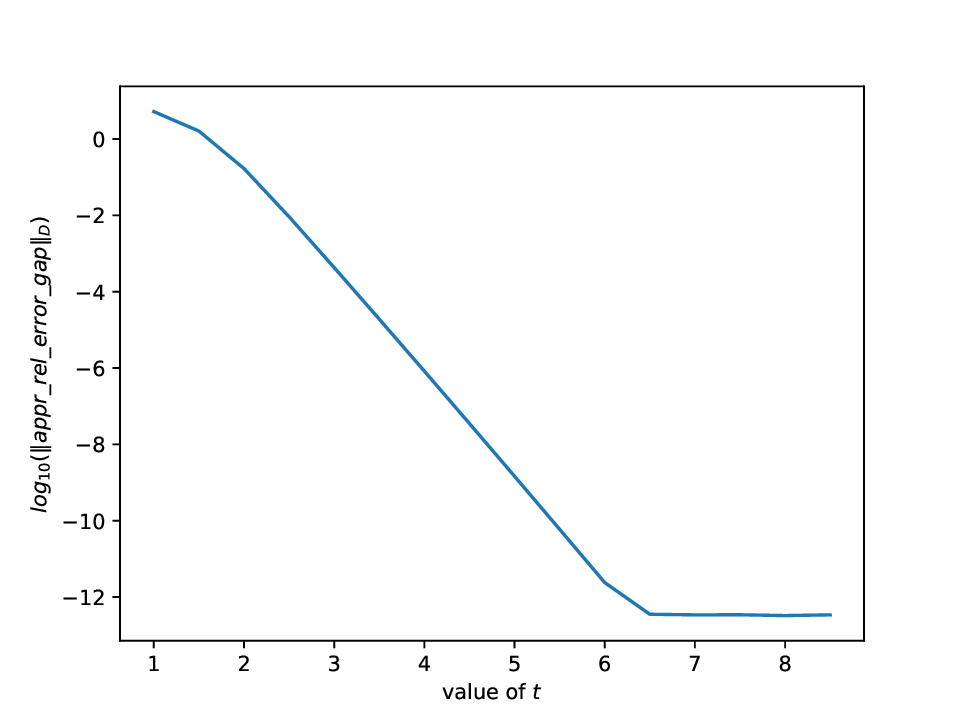}
\caption{Network of queues 2, approximate relative error gap versus $t$.}
\label{fig5}
\end{figure}

Figures \ref{fig::fig4} and \ref{fig5} are concerend with the setting of Section \ref{sec:network_queues2}.
Because the exact solution is unkown here, we modify our accuracy criteria.
Figure \ref{fig::fig4} shows the value of two accuracy measures for each state in  the truncation set $\{(i_1, i_2): i_1 \leq 50, i_2 \leq 50 \}$.
Figure \ref{fig4:figa} displays $\log_{10}$ of the approximate relative error gap which is defined as follows. 
\begin{align*}
appr\_rel\_error\_gap(x) = \frac{|upper\_b(x) - lower\_b(x)|}{|upper\_b(x)|}
\end{align*}
For the denominator, we did not use $lower\_b$ as it could change sign and become negative close to the boundary between $A$ and $A^c$.
Figure \ref{fig4:figb} displays $\log_{10}$ of the absolute error gap:
\begin{align*}
abs\_error\_gap(x) = |upper\_b(x) - lower\_b(x)|.
\end{align*}
Again, these figures suggest that close to the boundary between $A$ and $A^c$, the approximation error of $upper\_b, lower\_b$ is high.

Again, suppose we are interested in estimating the solution to Poisson's equation for a finite set $D  = \{(i_1,i_2) \in S: (i_1,i_2) \neq z, i_1\leq 20 , i_2\leq 20  \}$.
Figure \ref{fig5} displays $\|appr\_rel\_error\_gap\|_{D}$ as a function of $t$ where $t$ depicts the truncation set to be $A_t = \{(i_1,i_2) \in S: i_1\leq \lceil 20t \rceil, i_2\leq \lceil 20t \rceil \}$ and where $\|appr\_rel\_error\_gap\|_{D} = \max_{x \in D} appr\_rel\_error\_gap(x)$. Again, the figure indicates $\|appr\_rel\_error\_gap\|_{D}$ decreases as a function of $t$.

\bibliographystyle{apalike}
\bibliography{papers}

\end{document}